\documentclass[11pt]{article}

\usepackage{amsmath, amsfonts, amssymb,latexsym,wasysym,graphicx}
\input epsf.sty

\newtheorem{Lemma}{Lemma}[section]
\newtheorem{Proposition}[Lemma]{Proposition}

\newtheorem{Definition}[Lemma]{Definition}

\newcommand{\BEQ}{\begin{equation}}     
\newcommand{\BEA}{\begin{eqnarray}}
\newcommand{\BD}{\begin{displaymath}}
\newcommand{\EEQ}{\end{equation}}       
\newcommand{\EEA}{\end{eqnarray}}
\newcommand{\ED}{\end{displaymath}}
\newcommand{\al}{\alpha}                

\newcommand{\del}{\delta}

\newcommand{\eps}{\varepsilon}          

\newcommand{\R}{\mathbb{R}}
\newcommand{\C}{\mathbb{C}}

\newcommand{\Id}{{\mathrm{Id}}}

\def\proba{{\mathbb{P}}}
\def\esper{{\mathbb{E}}}

\def\Var{{\mathrm{Var}}}

\def\sgn{{\mathrm{sgn}}}

\newcommand{\eop}{\hfill $\Box$}        

\newcommand{\II}{{\rm i}}               
\renewcommand{\Re}{{\rm Re\ }}          
\renewcommand{\Im}{{\rm Im\ }}          
\newcommand{\half}{{1\over 2}}          
\renewcommand{\vec}[1]{\boldsymbol{#1}} 

\catcode`\@=11
\def\numberbysection{\@addtoreset{equation}{section}
        \def\theequation{\thesection.\arabic{equation}}}
\numberbysection

\begin{document}

\vspace*{1.5cm}
\begin{center}
{\Large \bf Moment estimates for solutions of linear stochastic differential
equations driven by analytic fractional Brownian motion
}
\end{center}

\vspace{2mm}
\begin{center}
{\bf  J\'er\'emie Unterberger}
\end{center}

\vspace{2mm}
\begin{quote}

\renewcommand{\baselinestretch}{1.0}
\footnotesize
{As a general rule, differential equations driven by a multi-dimensional irregular path $\Gamma$ are 
 solved  by constructing a rough path over $\Gamma$. The domain of definition -- and also
estimates -- of the solutions depend  on upper bounds for the rough path; these
general, deterministic estimates are too crude to  apply e.g. to   the solutions of stochastic
differential equations with linear coefficients
 driven by a Gaussian process with H\"older regularity $\alpha<1/2$.

We prove here (by showing convergence of Chen's series)
  that linear stochastic differential equations driven by analytic fractional Brownian motion
\cite{TinUnt,Unt08}
with arbitrary Hurst index $\alpha\in(0,1)$ may be solved on the closed upper half-plane, and that the solutions have finite variance.
 }
\end{quote}

\vspace{4mm}


\newpage


\section{Introduction}


 Assume
$\Gamma_t=(\Gamma_t(1),\ldots,\Gamma_t(d))$ is a smooth $d$-dimensional path, and $V_1,\ldots,V_d:\R^r\to\R^r$
be smooth vector fields. Then (by the classical Cauchy-Lipschitz theorem for instance) the differential
equation driven by $\Gamma$
\BEQ dy(t)=\sum_{i=1}^d V_i(y(t))d\Gamma_i(t) \label{eq:0:eq-dif} \EEQ
admits a unique solution with initial condition $y(0)=y_0$. The usual way to prove this is
by showing (by a functional  fixed-point theorem) that iterated integrals 
\BEQ y_n(t)\to
y_{n+1}(t):=y_0+\int_0^t \sum_i V_i(y_n(s))d\Gamma_i(s)  \label{eq:0:fixed-point} \EEQ
converge when $n\to \infty$.

Expanding this expression to all orders yields formally for an arbitrary analytic function $f$
\BEQ f(y_t)=f(y_s)+\sum_{n=1}^{\infty} \sum_{1\le i_1,\ldots,i_n\le d} \left[ V_{i_1}\ldots V_{i_n}f\right]
(y_s) {\bf \Gamma}^n_{ts}(i_1,\ldots,i_n), \label{eq:0:series} \EEQ
where 
\BEQ {\bf\Gamma}^n_{ts}(i_1,\ldots,i_n)=:=\int_s^t d\Gamma_{t_1}(i_1)\int_s^{t_1} d\Gamma_{t_2}(i_2)
\ldots \int_{s}^{t_{n-1}} d\Gamma_{t_n}(i_n),\EEQ
provided, of course, the series converges. By specializing to the identity function $f=\Id:\R^r\to\R^r$,
$x\to x$, one gets a series expansion for the solution $(y_t)$.

This formula, somewhat generalized, has been used in a variety of contexts:

\begin{enumerate}

\item Let 
\BEQ {\cal E}_V^{N,t,s}(y_s)=y_s+\sum_{n=1}^N  \sum_{1\le i_1,\ldots,i_n\le d} \left[ V_{i_1}\ldots V_{i_n}\Id\right]
(y_s) {\bf \Gamma}^n_{ts}(i_1,\ldots,i_n) \label{eq:0:truncated-series} \EEQ
be the $N$-th order truncation of (\ref{eq:0:series}). It may be interpreted as {\it one} iteration
of the numerical Euler scheme of order $N$, which is defined by 
\BEQ y_{t_k}^{Euler;D}:={\cal E}_V^{N,t_k,t_{k-1}}\circ\ldots\circ {\cal E}_V^{N,t_1,t_0}(y_0) \EEQ
for an arbitrary partition $D=\{0=t_0<\ldots<t_n=T\}$ of the interval $[0,T]$. When $\Gamma$
is only $\alpha$-H\"older with $\frac{1}{N+1}<\alpha\le \frac{1}{N}$, the iterated integrals 
${\bf \Gamma}^n(i_1,\ldots,i_n)$, $n=2,\ldots,N$ do not make sense a priori and must be substituted
with a {\it geometric rough path} over $\Gamma$. A {\it geometric rough path} over $\Gamma$ is a family
\BEQ \left( ({\bf \Gamma}^1_{ts}(i_1))_{1\le i_1\le d}, ({\bf\Gamma}^2_{ts}(i_1,i_2))_{1\le i_1,i_2\le d},\ldots,
({\bf\Gamma}^N_{ts}(i_1,\ldots,i_N)_{1\le i_1,\ldots,i_N\le d}) \right) \EEQ
of functions of two variables such that: 

\begin{itemize}
\item[(i)] $ {\bf\Gamma}^1_{ts}=\Gamma_t^1-\Gamma_s^1$;
\item[(ii)] ({\it{H\"older continuity}})
 each component of $\vec{\Gamma}^k$, $k=1,\ldots,N$
 is 
$k\alpha$-H\"older continuous, that is
to say, $\sup_{s\in\R} \left(\sup_{t\in\R} \frac{|{\bf \Gamma}_{ts}^k(i_1,\ldots,i_k)|}{|t-s|^{k\alpha}} \right)<\infty.$

\item[(iii)] ({\it{multiplicativity}}) letting $\del{\bf\Gamma}^k_{tus}:=
{\bf\Gamma}_{ts}^k-{\bf\Gamma}^k_{tu}-{\bf \Gamma}^k_{us}$, one requires
\BEQ
 \del\vec{\Gamma}^k_{tus}(i_1,\ldots,i_k) = \sum_{k_1+k_2=k} \vec{\Gamma}_{tu}^{k_1}(i_1,\ldots,i_{k_1}) \vec{\Gamma}_{us}^{k_2}(i_{k_1+1},\ldots,i_k). \label{eq:0:x}
\EEQ

\item[(iii)] ({\it{geometricity}}) 
\BEQ  {\bf \Gamma}^{n_1}_{ts}(i_1,\ldots,i_{n_1}) {\bf \Gamma}^{n_2}_{ts}(j_1,\ldots,j_{n_2}) 
 = 
\sum_{\vec{k}\in {\mathrm{Sh}}(\vec{i},\vec{j})} {\bf \Gamma}^{n_1+n_2}(k_1,\ldots,k_{n_1+n_2})  \label{eq:0:geo} \EEQ
where ${\mathrm{Sh}}(\vec{i},\vec{j})$ is the subset of permutations of $i_1,\ldots,i_{n_1},j_1,\ldots,j_{n_2}$
which do not change the orderings of $(i_1,\ldots,i_{n_1})$ and $(j_1,\ldots,j_{n_2})$.
\end{itemize}

Properties (i)-(iv) are true when $\Gamma$ is regular; the multiplicative
 property measures in some sense the
defect of additivity of iterated integrals, which is easy to measure when one represents these as geometric
quantities (areas, volumes, etc.)  Under these conditions, it is possible to integrate a $1$-form along the path $\Gamma$ (or, more
precisely, along the rough path $\bf\Gamma$); we refer the reader either to \cite{FV} or to \cite{Gu}.

It is also possible to solve differential equations driven by $\Gamma$ like (\ref{eq:0:eq-dif}), either
by using eq. (\ref{eq:0:fixed-point}) and a fixed-point theorem in a class of $\Gamma$-controlled processes
\cite{Gu}, or by
using the above Euler scheme \cite{FV}.  Assuming either (i) the vector fields $V$  and their derivatives up to order $N$ are bounded or
(ii) they are linear \footnote{Similar results are also expected when the vector fields $V$
are uniformly Lipschitz on $\R$ \cite{GubLej}. See also \cite{Lej} for investigations on the
same subject using a fixed point method.}, then the solution is globally defined, and the solution at time $T$ is bounded (i) by a polynomial
in $ |||{\bf\Gamma}|||$ or (ii) by
something like $\exp C |||{\bf\Gamma}|||^N$, where $|||{\bf\Gamma}|||=\max_{n=1,\ldots,N}
\sup_{1\le i_1,\ldots,i_n\le d}
\sup_{0\le s,t\le T} \frac{|{\bf\Gamma}^n_{ts}(i_1,\ldots,i_n)|}{|t-s|^{n\alpha}}$ is the
maximum of the H\"older norms. It seems there is no way (using either approach) to improve these
bounds in the general deterministic setting. Unfortunately, they do not give a control of the solution as a
stochastic process in the linear case (ii) when $\Gamma$ is a Gaussian process (such as fBm or analytic fBm, see below) with
H\"older regularity $\alpha<\half$, since $\esper \exp C(\Gamma_t-\Gamma_s)^2=\infty$ for $C$
large enough, and in any case $\esper \exp C|\Gamma_t-\Gamma_s|^N=\infty$ if $N\ge 3$.

\item Assume $\Gamma$ is a stochastic process, and let $P_t(f)=\esper f(y_t)$. When (\ref{eq:0:eq-dif})
is a diffusion driven by usual Brownian motion, $P_t$ is the associated semi-group operator. Assume
now $\Gamma$ is more general, for instance fBm or analytic fBm. Even
though the process in not Markov, the operator $P_t$ is interesting in itself.  The small-time
expansion of $P_t$ (corresponding to an arbitrary truncation of the above series)
 has been studied \cite{BauCou} when $\Gamma$ is fBm with Hurst index $\alpha>1/3$. When $ \alpha>1/2$, it
has been proved \cite{NNRT} that the series converges for functions and vector fields $V$ satisfying
somewhat drastic conditions.

\end{enumerate}

In any case, it seems difficult to get moment estimates for the solutions of stochastic differential
equations driven by stochastic processes $\Gamma$
 with H\"older regularity $\alpha<\half$. One of the reasons 
\cite{NNRT} is the difficulty of getting estimates for the iterated integrals ${\bf\Gamma}$; another reason
lies in the essence of the rough path method which relies on pathwise estimates; a third reason is, of course, that
the Chen series diverges even in the simplest cases (one-dimensional usual Brownian motion for instance) as soon
as the vector fields are unbounded and non-linear, e.g. quadratic.

\bigskip

In this article, we prove convergence of the series (\ref{eq:0:series}) when
the vector fields $V_i$ are linear and  $\Gamma$ is {\it analytic fBm} (afBm for short).
This process -- first defined in \cite{Unt08} --, depending on an index $\alpha\in(0,1)$, is a 
complex-valued process, a.s. $\kappa$-H\"older for every $\kappa<\alpha$, which has an analytic continuation to
the upper half-plane $\Pi^+:=\{z=x+\II y\ |\ x\in\R, y>0\}$.
 Its real part $(2\Re\Gamma_t, t\in\R)$ has the same law as fBm with Hurst index
$\alpha$. Trajectories of $\Gamma$ on the closed upper half-plane $\bar{\Pi}^+=\Pi^+\cup\R$
have the same regularity as those of fBm, namely, they are
$\kappa$-H\"older for every $\kappa<\alpha$. 
As shown in \cite{TinUnt}, the regularized rough path -- constructed by moving inside the
upper half-plane through an imaginary translation $t\to t+\II\eps$ -- converges in the limit $\eps\to 0$
to a geometric rough path over $\Gamma$ for {\it any} $\alpha\in(0,1)$, which makes it possible to produce
 strong, local pathwise solutions
of stochastic differential equations driven by $\Gamma$ with analytic coefficients.

We do not enquire about the convergence of the series (\ref{eq:0:series}) in the general case (as mentioned
before, it diverges e.g. when  $V$ is quadratic), but only in the
{\it linear} case. One obtains:

{\bf Main theorem.}

{\em Let $V_1,\ldots,V_d$ be linear vector fields on $\C^r$. Then the series (\ref{eq:0:series}) converges
in $L^2(\Omega)$ on the closed upper half-plane $\bar{\Pi}^+=\Pi^+\cup\R$. Furthermore, the solution
$(y_t)_{t\in\bar{\Pi}^+}$, defined as the limit of the series, has finite variance.
More precisely,  there exists a constant $C$ such that
\BEQ \esper |y_t-y_s|^{2}\le C |t-s|^{2\alpha},\quad s,t\in \bar{\Pi}^+. \label{eq:main-th} \EEQ

}

{\bf Notation.} Constants (possibly depending on $\alpha$) are generally denoted by $C,C',C_1,c_{\alpha}$
 and so on.


\section{Definition of afBm and first estimates}


We briefly recall to begin with the definition of the analytic fractional Brownian motion $\Gamma$, which
is a complex-valued process defined on the closed upper half-plane $\bar{\Pi}^+$ \cite{TinUnt}. 
Its introduction was initially motivated by the possibility to construct quite easily iterated integrals of $\Gamma$
by a contour deformation. Alternatively, its Fourier transform is supported on $\R_+$, which makes the regularization
procedure in \cite{Unt09a,Unt09b} void.

\begin{Proposition} \label{prop:1}

Let $\{\xi_k^1,\xi_k^2;\, k\ge 0\}$ be two families of independent standard Gaussian random variables,
 defined on a complete probability space $(\Omega,{\cal F},\proba)$,
 and for $k\ge 0$, set $\xi_k=\xi_k^1+\II\xi_k^2$.
 Consider the process $\Gamma'$ defined for $z\in\Pi^+$ by $\Gamma'_z=\sum_{k\ge 0}f_k(z)\xi_k$. Then:

\smallskip

\noindent
{\bf (1)}
$\Gamma'$ is a well-defined analytic process on $\Pi^+$, with Hermitian  covariance kernel
\BEQ \esper \Gamma'_z \Gamma'_w=0, \quad \esper \Gamma'_z\bar{\Gamma}'_w=\frac{\alpha(1-2\alpha)}{2\cos\pi\alpha} (-\II(z-\bar{w}))^{2\alpha-2}.\EEQ

\smallskip

\noindent
{\bf (2)}
Let $\gamma:(0,1)\to\Pi^+$ be any continuous path with endpoints $\gamma(0)=0$ and $\gamma(1)=z$, and set 
$\Gamma_z=\int_{\gamma}\Gamma'_u \, du$. Then $\Gamma$ is an analytic process on $\Pi^+$. Furthermore,
 as $z$ runs along any path in $\Pi^+$ going to $t\in\R$, the random variables $\Gamma_z$ converge 
almost surely to a random variable called again $\Gamma_t$.

\smallskip

\noindent
{\bf (3)}
The family $\{\Gamma_t;\, t\in\R\}$ defines a Gaussian centered complex-valued process, whose covariance function is given by:
$$ \esper[\Gamma_s \Gamma_t]=0, \quad 
\esper[\Gamma_s\bar\Gamma_t]=
\frac{e^{-\II\pi\al \,\sgn(s)}|s|^{2\al}
+e^{\II\pi\al \,\sgn(t)}|t|^{2\al}-e^{\II\pi\al \,\sgn(t-s)}|s-t|^{2\al}}{4\cos(\pi\al)}.
$$
The paths of this process are almost surely $\kappa$-H\"older for any $\kappa<\al$.

\smallskip

\noindent
{\bf (4)}
Both real and imaginary parts of $\{\Gamma_t;\, t\in\R\}$ are (non independent) fractional Brownian motions indexed by $\R$, with covariance given by
\begin{equation} \esper [\Re \Gamma_s \Im \Gamma_t]=-\frac{\tan\pi\alpha}{8} \left[ -\sgn(s) |s|^{2\alpha}+
\sgn(t) |t|^{2\alpha}-\sgn(t-s) |t-s|^{2\alpha} \right].\end{equation}

\end{Proposition}

\begin{Definition}

Let $Y_t:=\Re \Gamma_{it}$, $t\in\R_+$. More generally, $Y_t=(Y_t(1),\ldots,Y_t(d))$ is a vector-valued
process with $d$ independent, identically distributed components.

\end{Definition}

The above results imply that $Y_t$ is real-analytic on $\R_+^*$.

\begin{Lemma}

The infinitesimal covariance function of $Y_t$ is:
\BEQ \esper Y'_s Y'_t=\frac{\alpha(1-2\alpha)}{4\cos\pi\alpha} (s+t)^{2\alpha-2}.\EEQ

\end{Lemma}

{\bf Proof.} Let $X_t:=\Im\Gamma_{it}$. Since $\esper \Gamma_s\Gamma_t=0$, $(Y_s,s\ge 0)$ and $(X_s,s\ge 0)$
have same law, with covariance kernel $\esper Y_s Y_t=\esper X_s X_t=\half \Re \Gamma_{is}\bar{\Gamma}_{it}$.
Hence
\BEQ \esper [Y'_s Y'_t]=\half \Re\esper \Gamma'_{is}\bar{\Gamma}'_{it}=\frac{\alpha(1-2\alpha)}{4\cos\pi\alpha} (s+t)^{2\alpha-2}.\EEQ \hfill \eop

Note that $\esper Y'_s Y'_t>0$. From this simple remark follows (see proof of a similar
statement in \cite{NNRT} concerning usual 
fractional Brownian motion with Hurst index $\alpha>1/2$):

\begin{Lemma} \label{lem:Y}

Let ${\bf Y}^n_{ts}(i_1,\ldots,i_n)$, $n\ge 2$ be the iterated integrals of $Y$. Then there exists a constant $C>0$ such that
\BEQ \Var {\bf Y}^n_{ts}(i_1,\ldots,i_n) \le C \frac{(C|t-s|)^{2n\alpha}}{n!}.\EEQ

\end{Lemma}

{\bf Proof.} Let $\Pi$ be the set of all pairings $\vec{\pi}$ of the set $\{1,\ldots,2n\}$ such that
$\left( (k_1,k_2)\in\vec{\pi}\right)\Rightarrow \left(i_{k'_1}=i_{k'_2}\right)$,
 where $k'_1=k_1$ if $k_1\le n$, $k_1-n$ otherwise,
and similarly for $k'_2$. By Wick's formula,
\BEA   && \Var {\bf Y}^n_{ts}(i_1,\ldots,i_n) \nonumber\\
&& \qquad =\sum_{\vec{\pi}\in\Pi} \left(\int_s^t dx_1\ldots
\int_s^{x_{n-1}} dx_n\right)  \left(\int_s^t dx_{n+1}\ldots
\int_s^{x_{2n-1}}dx_{2n}\right) \nonumber\\
&& \qquad \qquad  \prod_{(k_1,k_2)\in\vec{\pi}} \esper[Y'_{x_{k_1}}Y'_{x_{k_2}}]. \label{eq:Wick1} \EEA
Since the process $Y'$ is positively correlated, and $\Pi$ is largest when all indices $i_1,\ldots,i_n$
are equal, one gets $\Var {\bf Y}^n_{ts}(i_1,\ldots,i_n)\le \Var {\bf Y}^n_{ts}(1,\ldots,1)$. On the other hand,
${\bf Y}^n_{ts}(1,\ldots,1)=\frac{1}{n!} (Y_t-Y_s)^n$, hence
\BEQ \Var{\bf Y}^n_{ts}(1,\ldots,1)=\frac{[\Var (Y_t-Y_s)]^n}{(n!)^2}\ .\ \frac{(2n)!}{2^n\ .\ n!} \le
\frac{[2\Var(Y_t-Y_s)]^n}{n!}. \EEQ

Now (assuming for instance $0<s<t$) 
\BEQ \Var (Y_t-Y_s)=c_{\alpha} \int_s^t \int_s^t (u+v)^{2\alpha-2} dudv\le c_{\alpha}
s^{2\alpha-2} (t-s)^2 \le c_{\alpha} (t-s)^{2\alpha} \EEQ
if $\frac{t}{2}\le s\le t$, and
\BEQ \Var (Y_t-Y_s)=\frac{c_{\alpha}}{2\alpha(2\alpha-1)} \left[ (2t)^{2\alpha}+(2s)^{2\alpha}-
2(t+s)^{2\alpha}\right]\le Ct^{2\alpha}\le C'(t-s)^{2\alpha} \EEQ
if $s<t/2$. Hence the result.  \hfill \eop


\section{Estimates for iterated integrals of $\Gamma$}


The main tool for the study of $\Gamma$ is the use of contour deformation. Iterated integrals of $\Gamma$
are particular cases of analytic iterated integrals, see \cite{Unt08} or \cite{TinUnt}. In particular,
the following holds:

\begin{Lemma}

Let $\gamma:(0,1)\to\Pi^+$ be the piecewise linear contour with affine parametrization defined by : 
\begin{itemize}
\item[(i)] $\gamma([0,1/3])=[s,s+i|\Re (t-s)|]$;
\item[(ii)] $\gamma([1/3,2/3])=[s+\II|\Re(t-s)|, t+\II |\Re(t-s)|]$;
\item[(iii)] $\gamma([2/3,1])=[t+\II|\Re(t-s)|,t]$.
\end{itemize}

If $z=\gamma(x)\in\gamma([0,1])$, we let $\gamma_z$ be the same path stopped at $z$, i.e.
$\gamma_z=\gamma([0,x])$, with the same parametrization. Then (letting $c_{\alpha}=\frac{\alpha(1-2\alpha)}{2\cos\pi\alpha}$)
\BEA &&  \Var {\bf \Gamma}^n_{ts}(i_1,\ldots,i_n)= \nonumber\\
 && \qquad c_{\alpha}^n
\sum_{\sigma\in\Sigma_I} \int_{\gamma} dz_1 \int_{\bar{\gamma}}d\bar{w}_1 (-\II(z_1-\bar{w}_{\sigma(1)}))^{2\alpha-2} \ .\ \int_{\gamma_{z_1}}dz_2\int_{\bar{\gamma}_{\bar{w}_1}} d\bar{w}_2 (-\II(z_2-\bar{w}_{\sigma(2)}))^{2\alpha-2} \ldots \nonumber\\
&& \qquad \qquad  \int_{\gamma_{z_{n-1}}}dz_n\int_{\bar{\gamma}_{\bar{w}_{n-1}}} d\bar{w}_n
(-\II(z_n-\bar{w}_{\sigma(n)}))^{2\alpha-2}  \label{eq:Gamma-gamma}  \EEA

where $\Sigma_I$ is the subset of permutations of $\{1,\ldots,n\}$ such that $(i_j=i_k)\Rightarrow
(\sigma(j)=\sigma(k))$.

\end{Lemma}

{\bf Proof.} Note first that, similarly to eq. (\ref{eq:Wick1}), 
\BEA &&\Var {\bf \Gamma}_{ts}^n(i_1,\ldots,i_n)=\sum_{\sigma\in\Sigma_I} \left(
\int_1^t dz_1\ldots \int_s^{z_{n-1}} dz_n\right) \left( \int_{\bar{s}}^{\bar{t}} d\bar{w}_1
\ldots \int_{\bar{s}}^{\bar{w}_{n-1}} d\bar{w}_n \right) \nonumber\\
  && \qquad \qquad \qquad \prod_{j=1}^n \esper \left[ \Gamma'_{z_j} \bar{\Gamma}'_{\bar{w}_{\sigma(j)}}
\right]
\EEA
(the difference with respect to eq. (\ref{eq:Wick1}) comes from the fact that contractions only operate between
$\Gamma$'s and $\bar{\Gamma}$'s, since $\esper [\Gamma_{z_j}\Gamma_{z_k}]=\esper [
\bar{\Gamma}_{\bar{w}_j}\bar{\Gamma}_{\bar{w}_k}]=0$ by Proposition \ref{prop:1}). 
Now the result comes from a deformation of contour, see \cite{Unt08}. \hfill \eop

\begin{Lemma} \label{lem:moment}

There exists a constant $C'$ such that, for every $s,t\in\bar{\Pi}^+=\Pi^+\cup\R$,

\BEQ \Var {\bf\Gamma}^n_{ts}(i_1,\ldots,i_n)\le \frac{(C'|t-s|)^{2n\alpha}}{n!}.\EEQ

\end{Lemma}

{\bf Proof.}  We assume (without loss of generality) that $\Im s\le \Im t$. If
$|\Im(t-s)|\ge c\Re |t-s|$  for some positive constant
$c$ (or equivalently $|\Re(t-s)|\le c' |t-s|$ for some $0\le c'<1$) then it is preferable to integrate along the straight line 
$[s,t]=\{z\in\C\ |\ z=(1-u)s+ut, \ 0\le u\le 1\}$ instead of $\gamma$, and use the parametrization
$y=\Im z$. If $z_1,z_2\in[s,t]$, $y_1=\Im z_1,y_2=\Im z_2$, then $|(-\II(z_1-\bar{z}_2))^{2\alpha-2}|
\le C(y_1+y_2)^{2\alpha-2}$, hence $\Var {\bf\Gamma}^n_{ts}(i_1,\ldots,i_n)\le C'^n \Var
 {\bf Y}^n_{y_2,y_1}
(i_1,\ldots,i_n)$, which yields the result by Lemma \ref{lem:Y}. So we shall assume that
$|\Re(t-s)|>c |t-s|$ for some constant $c>0$.

 Let us use as new variable the parametrization coordinate $x$ along $\gamma$. Then  formula 
(\ref{eq:Gamma-gamma}) reads
\BEA  && \Var {\bf \Gamma}^n_{ts}(i_1,\ldots,i_n)=c_{\alpha}^n \sum_{\sigma\in\Sigma_I} \int_0^1 dx_1 \int_0^1 dy_1 K'(x_1,y_{\sigma(1)})\ .\
\int_0^{x_1} dx_2 \int_0^{y_1} dy_2 K'(x_2,y_{\sigma(2)})\ldots \nonumber\\
&& \qquad \qquad  \int_0^{x_{n-1}}dx_n \int_0^{x_n} dy_n
K'(x_n,y_{\sigma(n)}),\EEA
where $K'(x,y)=(3|\Re(t-s)|)^{2} (3(x+y)|\Re(t-s)|+2\Im s)^{2\alpha-2}$ if $0<x,y<1/3$,
$(3|\Re(t-s)|)^{2} (3((1-x)+(1-y))|\Re(t-s)|+2\Im t)^{2\alpha-2}$ if $2/3<x,y<1$, and is bounded by
a constant times $|t-s|^{2\alpha}$ otherwise thanks to the condition $|\Re(t-s)|>c|t-s|$.
 Note that $(x+y)^{2\alpha-2}>2^{2\alpha-2}$ if
$0<x,y<1$. Hence (if $0<x,y<1$) $|K'(x,y)|\le (C_1|t-s|)^{2\alpha} \left[ (x+y)^{2\alpha-2} +((1-x)+(1-y))^{2\alpha-2}\right]$, which is (up to a coefficient) the infinitesimal covariance of $|t-s|^{\alpha}
(Y_x+\tilde{Y}_{1-x},0<x<1)$
if $\tilde{Y}\stackrel{(law)}{=} Y$ is independent of $Y$. A slight modification of the argument of Lemma
\ref{lem:Y} yields
\BEA  \Var {\bf \Gamma}^n_{ts}(i_1,\ldots,i_n) &\le&  (C_1|t-s|)^{2n\alpha}
\frac{[\Var (Y_1-Y_0)+\Var(\tilde{Y}_{1}-\tilde{Y}_{0})]^n}{(n!)^2} \ .\ \frac{(2n)!}{2^n\ .\ n!}  \nonumber\\
&\le& 
C_1^{2n\alpha} \ . \ \frac{(2C|t-s|)^{2n\alpha}}{n!}.\EEA
\hfill \eop

\section{Proof of main theorem}

We now prove the theorem stated in the introduction, which is really a simple corollary of Lemma
\ref{lem:moment}.

Let $C$ be the maximum of the matrix norms $|||V_i|||=\sup_{||x||_{\infty}=1} ||V_i x||_{\infty}$
for the supremum norm $||x||_{\infty}=\sup(|x_1|,\ldots,|x_r|)$.
 Rewrite eq. (\ref{eq:0:truncated-series}) as

\BEQ {\cal E}_V^{N,t,s}(y_s)=y_s+\sum_{n=1}^N  \sum_{1\le i_1,\ldots,i_n\le d}  a_{i_1,\ldots,i_n}
 {\bf \Gamma}^n_{ts}(i_1,\ldots,i_n).\EEQ

Then $||a_{i_1,\ldots,i_n}||_{\infty}\le C^n$ and $\esper |{\bf\Gamma}^n_{ts}(i_1,\ldots,i_n)|^2 \le 
\frac{(C|t-s|)^{2n\alpha}}{n!}$. Hence (by the Cauchy-Schwarz inequality)
\BEA \esper \left(  {\cal E}_V^{N,t,s}(y_s)-y_s \right)^2 &\le& \sum_{m,n=1}^N
\frac{(C''|t-s|)^{(m+n)\alpha}}{\sqrt{m!n!}} \nonumber\\
&=& \left( \sum_{m=1}^N \frac{(C''|t-s|)^{m\alpha}}{\sqrt{m!}} \right)^2 \le C'''. |t-s|^{2\alpha}\EEA
independently of $N$. The series obviously converges and yields eq. (\ref{eq:main-th}) for $p=1$.
\hfill \eop

It should be easy to prove along the same lines that the series defining $\esper |y_t-y_s|^{2p}$
converges for every $p\ge 1$,
 and that there exists a constant $C_p$ such that $\esper |y_t-y_s|^{2p}\le C_p |t-s|^{2\alpha p}$
for every $s,t\in\bar{\Pi}^+$.

\newpage


\end{document}